\documentclass{amsart}
\usepackage{bm,url}

\newcommand{\ra}{\rightarrow}
\newcommand{\C}{\mathbb{C}}
\newcommand{\R}{\mathbb{R}}
\newcommand{\Z}{\mathbb{Z}}
\newcommand{\Proj}{\mathbb{P}}
\newcommand{\bmu}{\bm{\mu}}
\renewcommand{\Re}{{\mathop{\rm Re}\nolimits}}
\newcommand{\Aut}{{\mathop{\rm Aut}\nolimits}}
\newcommand{\Sym}{{\mathop{\rm Sym}\nolimits}}
\newcommand{\co}{\colon}

\theoremstyle{plain}
\newtheorem{theorem}{Theorem}

\newtheorem{lemma}[theorem]{Lemma}
\newtheorem{proposition}[theorem]{Proposition}

\numberwithin{theorem}{section}

\theoremstyle{definition}

\numberwithin{equation}{section}

\title{The $D_4$ root system is not universally optimal}

\author{Henry Cohn}
\address{Microsoft Research\\
One Microsoft Way\\
Redmond, WA 98052-6399} \email{cohn@microsoft.com}

\author{John H.\ Conway}
\address{Department of Mathematics\\
Princeton University\\
Princeton, NJ 08544-1000} \email{conway@math.princeton.edu}

\author{Noam D.\ Elkies}
\address{Department of Mathematics\\
Harvard University\\
Cambridge, MA 02138} \email{elkies@math.harvard.edu}

\author{Abhinav Kumar}
\address{Department of Mathematics\\
Massachusetts Institute of Technology\\
Cambridge, MA 02139} \email{abhinav@math.mit.edu}

\thanks{Published in Experimental Mathematics \textbf{16} (2007), 313--320.}



\begin{document}

\begin{abstract}
We prove that the $D_4$ root system (equivalently, the set of
vertices of the regular 24-cell) is not a universally optimal
spherical code. We further conjecture that there is no universally
optimal spherical code of 24 points in $S^3$, based on numerical
computations suggesting that every 5-design consisting of
24 points in $S^3$ is in a 3-parameter family (which we
describe explicitly, based on a construction due to Sali) of
deformations of the $D_4$ root system.
\end{abstract}

\maketitle

\section{Introduction}\label{sec1}

In \cite{CK} the authors (building on work by Yudin, Kolushov,
and Andreev in \cite{Y,KY1,KY2,A1,A2}) introduce the notion of a
{\em universally optimal code} in $S^{n-1}$, the unit sphere in
$\R^n$. For a function $f\co [-1,1) \ra \R$ and a finite set $C
\subset S^{n-1}$, define the {\em energy} $E_f(C)$ by
$$
E_f(C) = \sum_{\genfrac{}{}{0pt}{}{c,c'\in C}{c\neq c'}}
f\big(\langle c, c'\rangle\big),
$$
where $\langle c, c'\rangle$ is the usual inner product.  We think
of $f$\/ as a potential function, and $E_f(C)$ as the potential
energy of the configuration $C$\/ of particles on $S^{n-1}$.  Note
that because each pair of points in $C$ is counted in both orders,
$E_f(C)$ is twice the potential energy from physics, but of course
this constant factor is unimportant.

A function $f \co [-1,1) \to \R$ is said to be {\em absolutely
monotonic} if it is smooth and it and all its derivatives are
nonnegative on $[-1,1)$.  A finite subset $C_0 \subset S^{n-1}$ is
said to be {\em universally optimal}\/ if $E_f(C_0) \leq E_f(C)$
for all $C\subset S^{n-1}$ with $\# C=\# C_0$ and all absolutely
monotonic $f$.  We say that $C_0$ is an {\em optimal spherical
code\/} if $t_{\max}(C_0) \leq t_{\max}(C)$ for all such $C$,
where
$$
t_{\max}(C) := \max_{\genfrac{}{}{0pt}{}{c,c'\in C}{c\neq c'}}
\langle c, c'\rangle
$$
is the cosine of the minimal distance of $C$. A universally
optimal code is automatically optimal (let $f(t)=(1-t)^{-N}$ or
$(1+t)^N$ for large $N$\/).

In \cite{CK}, linear programming bounds are applied to show that
many optimal codes are in fact universally optimal.  Notably
absent from this list is the $D_4$ root system $C_{D_4}$, which is
expected but not yet proved to be the unique optimal code of
size $24$ in $S^3$.  This root system can also be described as the
vertices of the regular \hbox{$24$-cell}.

It is shown in \cite{CK}
that the vertices of any regular polytope whose faces are
simplices form a universally optimal spherical code. The
dodecahedron, \hbox{$120$-cell}, and cubes in $\R^n$ with $n \ge
3$ are not even optimal spherical codes (see \cite{S}) and hence
cannot be universally optimal. Thus the \hbox{$24$-cell} was the
only remaining regular polytope.
Cohn and Kumar conjectured in an early draft of \cite{CK} that
$C_{D_4}$ was universally optimal, but reported the numerical
result that for the natural potential function $f(t) = (1-t)^{-1}$
there was another code $C\subset S^3$ with $\# C=24$ at which
$E_f$ has a local minimum only slightly larger than $E_f(C_{D_4})$
($668.1902+$, compared with $668$).

What makes this code
noteworthy is that in simulations of particle dynamics on $S^3$
under the potential function $f$ (along with a viscosity force to
remove kinetic energy and cause convergence to a local minimum for
$E_f$), $24$ particles converge more than 90\% of the time to $C$,
rather than to $C_{D_4}$. Similar effects appear to occur for
$f(t) = (1-t)^{-s}$ for other values of $s$.  In other words,
these codes have a much larger basin of attraction than $C_{D_4}$,
despite being suboptimal.

In this paper we give a simple description of a one-parameter
family of configurations $C_\theta$ that includes these codes, and
exhibit choices of $f$\/ (such as $f(t)=(1+t)^8$) and $\theta$ for
which $E_f(C_\theta) < E_f(C_{D_4})$. We thus disprove the
conjectured universal optimality of $C_{D_4}$.

We further conjecture that there is no universally optimal
spherical code of $24$ points in $S^3$.  Any such code would have
to be a \hbox{$5$-design}, because $C_{D_4}$ is.  Numerical
computations led us to a \hbox{$3$-parameter} family of such
designs that can be constructed using an approach introduced by
Sali in \cite{Sa}. The family contains $C_{D_4}$ as a special
case, and consists of deformations of $C_{D_4}$.

We exhibit these
designs and prove that, within the family, $C_{D_4}$ minimizes the
energy for every absolutely monotonic potential function, and is
the unique minimizer unless that function is a polynomial of
degree at most $5$. Our computations suggest that every
\hbox{$5$-design} of $24$ points in $S^3$ is in the new family. If
true, this would imply the nonexistence of a universally optimal
design of this size in $S^3$ because we already know that
$C_{D_4}$ is not universally optimal.

One way to think about the $D_4$ root system's lack of universal
optimality is that it explains how $D_4$ is worse than $E_8$.  The
$D_4$ and $E_8$ root systems are similar in many ways: they are
both beautiful, highly symmetrical configurations that seem to be
the unique optimal spherical codes of their sizes and dimensions.
However, one striking difference is that linear programming bounds
prove this optimality and uniqueness for $E_8$ but not for $D_4$
(see \cite{AB,BS,L,OS}).  This leads one to wonder what causes
that difference.  Is $D_4$ in some way worse than $E_8$?  Our
results in this paper show that the answer is yes: for $E_8$,
linear programming bounds prove universal optimality
(see \cite{CK}), while for $D_4$ universal optimality is not merely
unproved but in fact false.

\section{The codes $C_\theta$}\label{sec2}

We computed the $24 \times 24$ Gram matrix of inner products
between the points of the suboptimal but locally optimal
configuration mentioned above for the potential function
$f(t)=(1-t)^{-1}$. Each inner product occurred more than once,
suggesting that the configuration had some symmetry.  By studying
this pattern we eventually identified the configuration with a
code in the following family of $24$-point codes $C_\theta \subset
S^3$. (We are of course not the first to use this approach of
computing a code numerically and using its Gram matrix to detect
symmetries and then find good coordinates.  One recent case ---
also, as it happens, for codes in $S^3$ --- is \cite{SHC}, where
the method is called ``beautification.'')

We identify $\R^4$ with the complex vector space $\C^2$ so that
$$
S^3 = \{ (w_1,w_2) \in \C^2 : |w_1|^2 + |w_2|^2 = 1 \}.
$$
For $\theta \in \R / 2\pi\Z$ such that $\sin 2\theta \neq 0$ and
$\sin\theta \neq \cos\theta$ we set
$$
C_\theta := \{ (z,0), (0,w), (z\sin\theta, w\cos\theta),
(z\cos\theta, w\sin\theta) : z^3 = w^3 = 1 \}.
$$
Thus $C_\theta$ consists of $24$ unit vectors, namely $3+3$ of the
form $(z,0)$ or $(0,w)$ and $3^2+3^2$ of the form $(z\sin\theta,
w\cos\theta)$ or $(z\cos\theta, w\sin\theta)$.  Each of these
codes has $72$ symmetries (each complex coordinate can be
independently conjugated or multiplied by cube roots of unity, and
the two coordinates may be switched), forming a group $G$\/
isomorphic to the wreath product of the symmetric group $S_3$ with
$S_2$. This group does not act transitively: there are two orbits,
one consisting of the six points $(z,0)$ and $(0,w)$ and the other
consisting of the remaining $18$ points.

Listing all possible pairs $c,c' \in C_\theta$ with $c\neq c'$, we
find that there are in general $11$ possible inner products, with
multiplicities ranging from $18$ to $84$.  We thus compute that
\begin{align*}
E_f(C_\theta)  = {} &  18\bigl( f(0)+f(\sin 2\theta) \bigr)
+ 36 \bigl( f(\sin\theta) + f(\cos\theta) \bigr)
\\
& \phantom{} + 36\left( f\!\left(\sin^2\theta - \frac12
\cos^2\theta\right) + f\!\left(\cos^2\theta - \frac12
\sin^2\theta\right) \right)
\\
& \phantom{} + 72 \left( f\!\left(-\,\frac{\sin\theta}{2}\right)
+ f\!\left(-\,\frac{\cos\theta}{2}\right) + f\!\left(\frac{\sin
2\theta}{4}\right) + f\!\left(-\,\frac{\sin 2\theta}{2}\right)
\right)
\\
& \phantom{} + 84 f\!\left(-\frac12\right).
\end{align*}
For $C_{D_4}$ we have the simpler formula
$$
E_f(C_{D_4}) = 24 f(-1) + 192 \left(f\left(\frac12\right)+f\left(-\frac12\right)\right) + 144 f(0).
$$

\section{Failure of Universal Optimality}\label{sec3}

By  Theorem 9b in \cite[p.\ 154]{Wi}, an absolutely monotonic
function on $[-1,1)$ can be approximated, uniformly on compact
subsets, by nonnegative linear combinations of the absolutely
monotonic functions $f(t)=(1+t)^k$ with $k\in\{0,1,2,\dots\}$. To
test universal optimality of some spherical code $C_0$ it is thus
enough to test whether $E_f(C_0) \leq E_f(C)$ holds for all
$C\subset S^{n-1}$ with $\# C=\# C_0$ and each $f(t)=(1+t)^k$. We
wrote a computer program to compute $E_f(C_{D_4})$ and
$E_f(C_\theta)$, and plotted the difference $E_f(C_{D_4}) -
E_f(C_\theta)$ as a function of~$\theta$.

For $k \leq 2$ the plots suggested that $E_f(C_{D_4}) =
E_f(C_\theta)$ for all $\theta$.  This is easy to prove, either
directly from the formulas or more nicely by observing that
$C_{D_4}$ and $C_\theta$ are both spherical \hbox{$2$-designs}
(the latter because $G$\/ acts irreducibly on $\R^4$), so
$E_f(C_{D_4}) + 24 f(1)$ and $E_f(C_\theta) + 24 f(1)$ both equal
$24$ times the average of $c \mapsto f\big(\langle c,
c_0\rangle\big)$ over $S^3$ for any $c_0 \in S^3$.

For $k=3$, the plot suggested that $E_f(C_{D_4}) \leq
E_f(C_\theta)$, with equality at a unique value of $\theta$ in
$[0,\pi]$, numerically $\theta = 2.51674+$. We verified this
by using the rational parametrization
$$
\sin\theta = \frac{2u}{1+u^2}, \quad \cos\theta =
\frac{1-u^2}{1+u^2}
$$
of the unit circle, computing $E_f(C_{D_4}) - E_f(C_\theta)$
symbolically as a rational function of $u$, and factoring this
function. We found that
$$
E_f(C_{D_4}) - E_f(C_\theta) = -18 \frac{(u^6 - 6u^4 - 12u^3 +
3u^2 - 2)^2}{(u^2+1)^6},
$$
and thus that $E_f(C_{D_4}) \leq E_f(C_\theta)$, with equality if
and only if $u$ is a root of the sextic $u^6 - 6u^4 - 12u^3 + 3u^2
- 2$. This sextic has two real roots,
$$
u = -(0.51171+), \quad u=3.09594-,
$$
which yield the two permutations of $\{\sin\theta,\cos\theta\} =
\{ 0.58498-, -(0.81105+) \}$ and thus give rise to a unique code
$C_\theta$ with $E_f(C_{D_4}) = E_f(C_\theta)$. This code is
characterized more simply by the condition that $\sin\theta +
\cos\theta$ is a root of the cubic $3y^3-9y-2=0$, or better yet
that $\sin^3\theta+\cos^3\theta=-\frac13$. The latter formulation also
lets us show that this is the unique $C_\theta$ that is a
spherical \hbox{$3$-design}: the cubics on $\R^4$ invariant
under $G$\/ are the multiples of $\Re(w_1^3) + \Re(w_2^3)$, and
the sum of this cubic over $C_\theta$ is
$6+18(\sin^3\theta+\cos^3\theta)$. Since $C_{D_4}$ and this
particular $C_\theta$ are both \hbox{$3$-designs}, they
automatically minimize the energy for any potential function that
is a polynomial of degree at most $3$.  We must thus try $k>3$ if
we are to show that $C_{D_4}$ is not universally optimal.

For $k=4$ through $k=7$ the plot indicated that $E_f(C_\theta)$
comes near $E_f(C_{D_4})$ for $\theta \approx 2.52$ but stays
safely above $E_f(C_{D_4})$ for all $\theta$, which is easily
proved using the rational parametrization. (For $k=4$ and $k=5$ we
could also have seen that $C_{D_4}$ minimizes $E_f$ by noting that
$C_{D_4}$ is a \hbox{$5$-design}.) However, for $k=8$ the minimum
value of $E_f(C_\theta)$, occurring at $\theta = 2.529367746+$, is
$5064.9533+$, slightly but clearly smaller than $E_f(C_{D_4}) =
5065.5$. That is, this $C_\theta$ is a better code than $C_{D_4}$
for the potential function $(1+t)^8$, so $C_{D_4}$ is not optimal
for this potential function and hence not universally optimal.

The maximum value of $E_f(C_{D_4}) - E_f(C_\theta)$ for
$f(t)=(1+t)^k$ remains positive for $k=9,10,11,12,13$, attained at
values of $\theta$ that slowly increase from $\theta = 2.52937-$
for $k=8$ to $\theta = 2.54122-$ for $k=13$. Each of these is
itself enough to disprove the conjecture that $C_{D_4}$ is
universally optimal.  (Another natural counterexample is
$f(t)=e^{6t}$ with $\theta = 2.53719+$.)

We found no further solutions of $E_f(C_{D_4}) > E_f(C_\theta)$
with $k > 13$.  It is clear that $E_f(C_{D_4}) < E_f(C_\theta)$
must hold for all $\theta$ if $k$\/ is large enough, because
$t_{\max}(C_\theta) > t_{\max}(C_{D_4}) = \frac12$: the smallest value
$t_0$ of $t_{\max}(C_\theta)$ is $(\sqrt{7}-1)/3 = 0.54858+$,
occurring when either
$$
t_0 = \sin\theta = \cos^2\theta - \frac12\sin^2\theta
$$ with
$\theta = 2.56092+$ or
$$
t_0 = \cos\theta = \sin^2\theta - \frac12\cos^2\theta
$$
with $\theta = 5.29305+$. Quantifying what ``large enough'' means,
and combining the resulting bound with our computations for
smaller $k$, we obtained the following result:

\begin{proposition} \label{prop:summary}
For $8 \le k \le 13$, there exists a choice of $\theta$ for which
$$
E_f(C_\theta) < E_f(C_{D_4})
$$
when $f(t) = (1+t)^k$.  For other nonnegative integers $k$, no such
$\theta$ exists.
\end{proposition}

\begin{proof}
If $k$ is large enough that
$$
18 f\!\left(\frac{\sqrt{7}-1}{3}\right) > 24 f(-1) + 192
\left(f\!\left(\frac12\right)+f\!\left(-\frac12\right)\right) + 144 f(0),
$$
then $E_f(C_{D_4}) < E_f(C_\theta)$. This criterion is wasteful,
but we have no need of sharper inequalities.  Calculation shows
that this criterion holds for all $k \ge 75$. This leaves only
finitely many values of $k$. For each of them, we use the rational
parametrization of the unit circle to translate the statement of
Proposition \ref{prop:summary} into the assertion of existence or
nonexistence of a real solution of a polynomial in $\Z[u]$, which
can be confirmed algorithmically using Sturm's theorem. Doing so
in each case completes the proof.
\end{proof}

We cannot rule out the possibility that there exists a universally
optimal \hbox{$24$-point} code in $S^3$, but it seems exceedingly
unlikely.  If $C_{D_4}$ is the unique optimal spherical code, as
is widely believed, then no universally optimal code can exist.
The same conclusion follows from the conjecture in the next
section.

It is still natural to ask which configuration minimizes each
absolutely monotonic potential function.  We are unaware of any
case in which another code beats $C_{D_4}$ and all the codes
$C_\theta$ for some absolutely monotonic potential function, but
given the subtlety of this area we are not in a position to make
conjectures confidently.

\section{New Spherical 5-Designs}\label{sec4}

Spherical designs are an important source of minimal-energy
configurations: a spherical $\tau$-design automatically minimizes
the potential energy for $f(t) = (1+t)^k$ with $k \le \tau$.
Conversely, if an $N$-point spherical $\tau$-design exists in
$S^{n-1}$, then every $N$-point configuration in $S^{n-1}$ that
minimizes the potential function $f(t) = (1+t)^\tau$ must be a
$\tau$-design.  Thus, when searching for universally optimal
configurations, it is important to study $\tau$-designs with $\tau$
as large as possible.

For $24$ points in $S^3$, the $D_4$ root system forms a
$5$-design. By Theorem 5.11 in \cite{DGS}, every $6$-design must
have at least $30$ points, so $24$ points cannot form a
\hbox{$6$-design}. Counting degrees of freedom suggests that
$24$-point $4$-designs are plentiful, but $5$-designs exist only
for subtler reasons.  One can search for them by having a computer
minimize potential energy for $f(t) = (1+t)^5$.  Here, we report
on a three-dimensional family of \hbox{$5$-designs} found by this
method. The $D_4$ root system is contained in this family, and all
the designs in the family can be viewed as deformations of
$C_{D_4}$. We conjecture that there are no other \hbox{$24$-point}
spherical \hbox{$5$-designs} in $S^3$. We shall show that this
conjecture implies the nonexistence of a universally optimal
\hbox{$24$-point} code in $S^3$.

Our construction of \hbox{$5$-designs} slightly generalizes a
construction of Sloane, Hardin, and Cara for the \hbox{$24$-cell}
(Construction 1 and Theorem 1 in \cite{SHC}). The
Sloane--Hardin--Cara construction also works for certain other
dimensions and numbers of points, and can be further generalized
using our more abstract approach. For example, one can construct a
family of designs in $S^{2n-1}$ from a design in $\C\Proj^{n-1}$.
We plan to treat further applications in a future paper.  See also
the final paragraph of this section.

Fix an ``Eisenstein structure'' on the $D_4$ root lattice, that
is, an action of $\Z[r]$, where $r=e^{2\pi i/3}$ is a cube root of
unity. (It is enough to specify the action of $r$, which can be
any element of order $3$ in $\Aut(D_4)$ that acts on $\R^4$ with
no nonzero fixed points; such elements constitute a single
conjugacy class in $\Aut(D_4)$.) Then $\R^4$ is identified with
$\C^2$, with the inner product given by
$$
\big\langle(z_1,z_2),(\zeta_1,\zeta_2) \big\rangle =
\Re\big(z_1\overline{\zeta}_1 + z_2 \overline{\zeta}_2\big).
$$
The group $\bmu_6$ of sixth roots of unity (generated by $-1$ and
$r$) acts on $C_{D_4}$ and partitions its $24$ points into four
hexagons centered at the origin; call them $H_0, H_1, H_2, H_3$.
In coordinates, we may take
\begin{align*}
H_0 & = \{(w,0) : w \in \bmu_6\}, \\
H_1 & =  \{(uiw,tiw) : w \in \bmu_6\}, \\
H_2 & =  \{(uiw,rtiw) : w \in \bmu_6\},  \\
H_3 & =  \{(uiw, \overline{r}tiw) : w \in \bmu_6\},
\end{align*}
where $u = \sqrt{1/3}$ and $t = \sqrt{2/3}$. For any complex
numbers $a_0,a_1,a_2,a_3$ with $|a_0|=|a_1|=|a_2|=|a_3|=1$, we
define
$$
D(a_0,a_1,a_2,a_3) = a_0 H_0 \cup a_1 H_1 \cup a_2 H_2 \cup a_3
H_3.
$$
We claim that $D(a_0,a_1,a_2,a_3)$ is a $5$-design.  That is, we
claim that for every polynomial $P$\/ of degree at most $5$
on $\R^4$, the average of $P(c)$ over $c \in D(a_0,a_1,a_2,a_3)$
equals the average of $P(c)$ over $c \in S^3$. It is sufficient to
prove that the average is independent of the choice of
$a_0,a_1,a_2,a_3$, because $D(1,1,1,1) = C_{D_4}$ is already known
to be a \hbox{$5$-design}.  But this is easy: for each $m \in \{
0,1,2,3\}$, the restriction of $P$\/ to the plane spanned by $H_m$
is again a polynomial of degree at most $5$, and $a_m H_m$ is a
{$5$-design} in the unit circle of this plane, so the average of
$P$\/ over $a_m H_m$ is the average of $P$\/ over this unit
circle, independent of the choice of $a_m$.

This construction via rotating hexagons is a special case of
Lemma 2.3 in \cite{Sa}, where that idea is applied to prove that
many spherical designs are not rigid.  Sali rotates a single
hexagon to prove that $C_{D_4}$ is not rigid, but he does not
attempt a complete classification of the $24$-point $5$-designs.

It is far from obvious that there is no other way to perturb the
$D_4$ root system to form a $5$-design. For example, if there were
two disjoint hexagons in $D_4$ that did not come from the same
choice of Eisenstein structure as above, then rotating them
independently would produce $5$-designs not in our family.
However, one can check via a counting argument that every pair of
disjoint hexagons does indeed come from some common Eisenstein
structure. This supports our conjecture that there are no other
$24$-point spherical $5$-designs in $S^3$.

The family of \hbox{$5$-designs} of the form $D(a_0,a_1,a_2,a_3)$
is three-dimensional, for the following reason. One of the four
parameters $a_0,a_1,a_2,a_3$ is redundant, because for every
$\alpha \in \C^*$ with $|\alpha|=1$ we have $D(\alpha a_0, \alpha
a_1, \alpha a_2, \alpha a_3) \cong D(a_0,a_1,a_2,a_3)$. We may
thus assume $\alpha_0=1$. We claim that for each $(a_1,a_2,a_3)$
there are only finitely many $(a'_1,a'_2,a'_3)$ such that
$D(1,a_1,a_2,a_3) \cong D(1,a'_1,a'_2,a'_3)$.  If this were not
true, there would be an infinite set of designs $D(1,b_1,b_2,b_3)$
equivalent under automorphisms of $S^3$ that stabilize $H_0$
pointwise. But this is impossible, because such an automorphism
must act trivially on the first coordinate $z_1$.  Hence our
\hbox{$5$-designs} constitute a three-dimensional family, as
claimed.

Some other known spherical designs can be similarly generalized.
For instance, the $7$-design of $48$ points in $S^3$, obtained
in \cite{SHC} from two copies of $C_{D_4}$, has a decomposition
into six regular octagons, which can be rotated independently to
yield a five-dimensional family of $7$-designs.

It is also fruitful to take a more abstract approach.  A
$24$-point design with a $\bmu_6$ action is specified by four
points, one in each orbit.  Our new designs are characterized by
the condition that under the natural map $\C^2 \setminus \{(0,0)\}
\to \C\Proj^1$ given by $(z_1,z_2) \mapsto z_1/z_2$, the four
points must map to the vertices of a regular tetrahedron (if we
identify $\C\Proj^1$ with $S^2$ via stereographic projection, with
$\C\Proj^1 = \C \cup \{\infty\}$ and $S^2$ a unit sphere centered
at the origin).  In slightly different language, we have specified
the image of the design under the Hopf map $S^3 \to S^2$.  The
fact that the regular tetrahedron is a spherical $2$-design in
$S^2$ plays a crucial role, and can be used to prove that this
construction yields $5$-designs.  Likewise, the six octagons that
make up each of our \hbox{$7$-designs} map to the vertices of a
regular octahedron, which is a spherical {$3$-design}. Again, we
intend to discuss this approach in more detail in a future paper.

\section{Optimality of $C_{D_4}$ among new $5$-designs}\label{sec5}

In this section we prove that among all these new $5$-designs, the
$24$-cell minimizes potential energy for each absolutely monotonic
potential function. As before, it is sufficient to do this for
$f(t)=(1+t)^k$. For $k \le 5$ this follows immediately from the
spherical design property, and for $k>5$ we will show directly
that $C_{D_4}$ is the unique minimizer.

Within each hexagon $H_m$, the six points are in the same relative
position in each design, and thus make the same contribution to
the potential energy.  Hence it suffices to show that the
potential energy between each pair of hexagons is separately
minimized for the $D_4$ configuration.

Let $a_0=1$, $a_1 = ie^{i\theta}$, $a_2 = ie^{i\phi}$, and
$a_3=ie^{i\psi}$. Because of the sixfold rotational symmetry of
each $H_m$, the angles $\theta$, $\phi$, and $\psi$ are determined
only modulo $\pi/3$. In particular, $\theta=\phi=\psi=\pi/6$
yields the $24$-cell (because $\pi/2 \equiv \pi/6 \pmod{\pi/3}$).

First, consider the pair $(H_0,H_1)$. We find that the inner
products between the points of $H_0$ and $a_1H_1$ are $(1/\sqrt{3})
\cos(\theta + j\pi/3)$ with $0 \leq j \leq 5$, each repeated
six times. By Lemma \ref{lemma:min} below, the sum is minimized
exactly when $\theta \equiv \pi/6  \pmod{\pi/3}$.  Similarly,
considering $(H_0,H_2)$ and $(H_0,H_3)$ shows that the
corresponding contributions to potential energy are minimized when
$\phi \equiv \pi/6 \pmod{\pi/3}$ or $\psi \equiv \pi/6
\pmod{\pi/3}$, respectively.

Next, consider the pair $(H_1,H_2)$. The dot products possible are
of the form
$$
\Re\big((u^2+\overline{r}t^2)a_1\overline{a}_2w^j\big) =
\frac{1}{\sqrt{3}}\cos\left(\frac{3\pi}{2}+\theta-\phi + \frac{j\pi}{3}\right)
$$
with $0 \le j \le 5$, each repeated six times. Once again we
conclude from Lemma \ref{lemma:min} that the potential energy is
minimized when $\theta \equiv \phi \pmod{\pi/3}$.  Similarly,
considering the remaining pairs $(H_1,H_3)$ and $(H_2,H_3)$ shows
that $\theta \equiv \phi \equiv \psi \pmod{\pi/3}$. Thus, it will
follow from the lemma below that for each $k \ge 6$, the $D_4$
configuration with $\theta=\phi=\psi=\pi/6$ is the unique code in
this family that minimizes the potential energy under the
potential function $f(t)=(1+t)^k$.

\begin{lemma} \label{lemma:min}
Let $k$ be a nonnegative integer. When $k \ge 6$, the function
$$
\theta \mapsto \sum_{j=0}^5
\left(1+\frac{\cos(\theta+j\pi/3)}{\sqrt{3}}\right)^k
$$
has a unique global minimum within $[0,\pi/3]$, which occurs at
$\theta= \pi/6$.  When $k \le 5$, the function is constant.
\end{lemma}

\begin{proof}
We must show that the coefficient of $y^k$ in the generating
function
$$
\sum_{j=0}^5\sum_{k=0}^\infty
\left(\left(1+\frac{\cos(\theta+j\pi/3)}{\sqrt{3}}\right)^k -
\left(1+\frac{\cos(\pi/6+j\pi/3)}{\sqrt{3}}\right)^k \right)y^k
$$
is zero if $k \le 5$ or $\theta \equiv \pi/6 \pmod{\pi/3}$ and
strictly positive otherwise.  Explicit computation using the sum
of a geometric series shows that the generating function equals
$$
\frac{y^6\big(\!\cos^2\theta\big)\big(4\cos^2\theta-3\big)^2}{216}
\left(\frac{1}{1-y} + \frac{2}{2-y} +
\frac{2}{2-3y}\right)\prod_{j=0}^5
\frac{1}{1-y\left(1+\frac{\cos(\theta+j\pi/3)}{\sqrt{3}}\right)}.
$$
The factor of
$\big(\!\cos^2\theta\big)\big(4\cos^2\theta-3\big)^2$ vanishes iff
$\theta \equiv \pi/6 \pmod{\pi/3}$ and is positive otherwise.
Clearly the factor
$$
\frac{1}{1-y} + \frac{2}{2-y} + \frac{2}{2-3y}
$$
has positive coefficients, as does
$$
\prod_{j=0}^5
\frac{1}{1-y\left(1+\frac{\cos(\theta+j\pi/3)}{\sqrt{3}}\right)},
$$
because $1+\cos(\theta+j\pi/3)/\sqrt{3} >0$ for all $j$.  It
follows that their product has positive coefficients, and taking
the factor of $y^6$ into account completes the proof.
\end{proof}

\section{Local Optimality}\label{sec6}

So far, we have not addressed the question of whether our new codes
are actually local minima for energy.  Of course that is not
needed for our main result, because they improve on the $24$-cell
regardless of whether they are locally optimal, but it is an
interesting question in its own right.

For the codes $C_\theta$ this question appears subtle, and we do
not resolve it completely.  To see the issues involved, consider
the case of $f(t)=(1-t)^{-1}$.  As $\theta$ varies, the lowest
energy obtained is $668.1920+$ when $\theta = 2.5371+$.  That code
appears to be locally minimal among all codes, based on
diagonalizing the Hessian matrix numerically, but we have not
proved it. By contrast, the other two local minima within the
family $C_\theta$ (with energy $721.7796+$ at $\theta =
-(2.0231+)$ and energy $926.3218+$ at $\theta = 0.5320+$) are
critical points but definitely not local minima among all codes;
the Hessians have $22$ and $36$ negative eigenvalues,
respectively.

We do not know a simple criterion that predicts whether a local
minimum among the codes $C_\theta$ as $\theta$ varies will prove
to be a local minimum among all codes, but it is not hard to prove
that every critical point in the restricted setting is also an
unrestricted critical point.  Specifically, a short calculation
shows that for every code $C_\theta$ and every smooth potential
function, the gradient of potential energy on the space of all
configurations lies in the tangent space of the subspace
consisting of all the codes $C_\theta$.  It follows immediately
that if the derivative with respect to $\theta$ of potential
energy vanishes, then the gradient vanishes as well.  Furthermore,
such critical points always exist: starting at an arbitrary code
$C_\theta$ and performing gradient descent will never leave the
space of such codes and will always end at a critical point.

At this point one may wonder whether it is even clear that the
regular $24$-cell is a local minimum for all absolutely monotonic
potential functions.  It is straightforward to show that it is a
critical point, but we know of no simple proof that it is actually
a local minimum.  The best proof we have found is the following
calculation.

For each of the $24$ points, choose an orthonormal basis of the
tangent space to $S^3$ at that point, and compute the Hessian
matrix of potential energy with respect to these coordinates. Its
eigenvalues depend on the potential function, but the
corresponding eigenspaces do not.  There is a simple reason for
that, although we will not require this machinery. Consider the
space $\Sym^{24}(S^3)$ of all unordered sets of $24$ points in
$S^3$. The symmetry group of the $24$-cell acts on the tangent
space to $\Sym^{24}(S^3)$ at the point corresponding to the
$24$-cell, and this representation breaks up as a direct sum of
irreducible representations.  On each nontrivial irreducible
representation the Hessian has a single eigenvalue, and these
subspaces do not depend on the potential function. In practice,
the simplest way to calculate the eigenspaces is not to use
representation theory, but rather to find them for one potential
function and then verify that they are always eigenspaces.

If the potential function is $f \co [-1,1) \to \R$, then the
eigenvalues of the Hessian are
\begin{align*}
&0,\\
&2f''\!\left(\frac{1}{2}\right)+8f''(0)+2f''\!\left(-\frac{1}{2}\right)-12f'\!\left(\frac{1}{2}\right)+12f'\!\left(-\frac{1}{2}\right),\\
&2f''\!\left(\frac{1}{2}\right)+4f''(0)+6f''\!\left(-\frac{1}{2}\right)-8f'\!\left(\frac{1}{2}\right)-4f'(0)+8f'\!\left(-\frac{1}{2}\right)+4f'(-1),\\
&5f''\!\left(\frac{1}{2}\right)+4f''(0)+3f''\!\left(-\frac{1}{2}\right)-14f'\!\left(\frac{1}{2}\right)+8f'(0)+2f'\!\left(-\frac{1}{2}\right)+4f'(-1),\\
&6f''\!\left(\frac{1}{2}\right)+6f''\!\left(-\frac{1}{2}\right)-12f'\!\left(\frac{1}{2}\right)+12f'\!\left(-\frac{1}{2}\right),\\
&2f''\!\left(\frac{1}{2}\right)+4f''(0)+6f''\!\left(-\frac{1}{2}\right)+4f'\!\left(\frac{1}{2}\right)+8f'(0)+20f'\!\left(-\frac{1}{2}\right)+4f'(-1),\\
&6f''\!\left(\frac{1}{2}\right)+8f''(0)+6f''\!\left(-\frac{1}{2}\right)-4f'\!\left(\frac{1}{2}\right)+4f'\!\left(-\frac{1}{2}\right), \\
&8f''\!\left(\frac{1}{2}\right)+4f''(0)-8f'(1/2)-4f'(0)+8f'\!\left(-\frac{1}{2}\right)+4f'(-1),
\end{align*}
with multiplicities $6$, $9$, $16$, $8$, $12$, $4$, $9$, and $8$,
respectively.

One mild subtlety is that $0$ is always an eigenvalue, so one
might worry that the second derivative test is inconclusive.
However, note that the potential energy is invariant under the
action of the $6$-dimensional Lie group $\operatorname{O}(4)$, which yields the
$6$ eigenvalues of $0$.  In such a case, if all other eigenvalues
are positive, then local minimality still holds, for the following
reason. Notice that $\operatorname{O}(4)$ acts freely on the space of ordered
$24$-tuples of points in $S^3$ that span $\R^4$, and it acts
properly since $\operatorname{O}(4)$ is compact. The quotient space is therefore
a smooth manifold, and the positivity of the remaining eigenvalues
suffices for the potential energy to have a strict local minimum
on the quotient space.

To complete the proof, we need only consider $f(t) = (1+t)^k$ with
$k \in \{0,1,2,\dotsc\}$. For $k \le 5$ the other eigenvalues are
not all positive (some vanish), but because the $24$-cell is a
spherical $5$-design it is automatically a global minimum for
these energies.  For $k \ge 6$ one can check that all the other
eigenvalues are positive. That is obvious asymptotically, because
they grow exponentially as functions of $k$; to prove it for all
$k \ge 6$ one reduces the problem to a finite number of cases and
checks each of them.  It follows that the regular $24$-cell
locally minimizes potential energy for each absolutely monotonic
potential function, and it is furthermore a strict local minimum
(modulo orthogonal transformations) unless the potential function
is a polynomial of degree at most $5$.

\section*{Acknowledgments}

This work was begun during a workshop at the American Institute of
Mathematics. The numerical and symbolic computations reported here
were carried out with PARI/GP (see \cite{PARI}) and Maple.

We are grateful to Eiichi Bannai for providing helpful feedback
and suggestions.

NDE was supported in part by the National Science Foundation
and AK was supported by a summer internship in the theory group at
Microsoft Research and a Putnam fellowship at Harvard University.

\end{document}